# A method of moments estimator of tail dependence

JOHN H.J. EINMAHL[1,*] , ANDREA KRAJINA[1,**] and JOHAN SEGERS[2]

[1] *Department of Econometrics & OR and CentER, Tilburg University, P.O. Box 90153, 5000 LE Tilburg, The Netherlands. E-mail:* *j.h.j.einmahl@uvt.nl, **a.krajina@uvt.nl*
[2] *Institut de statistique, Université catholique de Louvain, Voie du Roman Pays, 20, B-1348 Louvain-la-Neuve, Belgium. E-mail: johan.segers@uclouvain.be*

In the world of multivariate extremes, estimation of the dependence structure still presents a challenge and an interesting problem. A procedure for the bivariate case is presented that opens the road to a similar way of handling the problem in a truly multivariate setting. We consider a semi-parametric model in which the stable tail dependence function is parametrically modeled. Given a random sample from a bivariate distribution function, the problem is to estimate the unknown parameter. A method of moments estimator is proposed where a certain integral of a nonparametric, rank-based estimator of the stable tail dependence function is matched with the corresponding parametric version. Under very weak conditions, the estimator is shown to be consistent and asymptotically normal. Moreover, a comparison between the parametric and nonparametric estimators leads to a goodness-of-fit test for the semiparametric model. The performance of the estimator is illustrated for a discrete spectral measure that arises in a factor-type model and for which likelihood-based methods break down. A second example is that of a family of stable tail dependence functions of certain meta-elliptical distributions.

## 1. Introduction

A bivariate distribution function $F$ with continuous marginal distribution functions $F_1$ and $F_2$ is said to have a *stable tail dependence function* $l$ if for all $x \geq 0$ and $y \geq 0$, the following limit exists:

$$\lim_{t \to 0} t^{-1} \mathbb{P}\{1 - F_1(X) \leq tx \text{ or } 1 - F_2(Y) \leq ty\} = l(x, y); \tag{1.1}$$

see [6, 15]. Here, $(X, Y)$ is a bivariate random vector with distribution $F$.

The relevance of condition (1.1) comes from multivariate extreme value theory: if $F_1$ and $F_2$ are in the max-domains of attraction of extreme value distributions $G_1$ and $G_2$ and if (1.1) holds, then $F$ is in the max-domain of attraction of an extreme value distribution $G$ with marginals $G_1$ and $G_2$ and with copula determined by $l$; see Section 2 for more details.

Inference problems on multivariate extremes therefore generally separate into two parts. The first one concerns the marginal distributions and is simplified by the fact that univariate extreme value distributions constitute a parametric family. The second one concerns the dependence structure in the tail of $F$ and forms the subject of this paper. In particular, we are interested in the estimation of the function $l$. The marginals will not be assumed to be known and will be estimated nonparametrically. As a consequence, the new inference procedures are rank-based and therefore invariant with respect to the marginal distribution, in accordance with (1.1).

The class of stable tail dependence functions does not constitute a finite-dimensional family. This is an argument for nonparametric, model-free approaches. However, the accuracy of these nonparametric approaches is often poor in higher dimensions. Moreover, stable tail dependence functions satisfy a number of shape constraints (bounds, homogeneity, convexity; see Section 2) which are typically not satisfied by nonparametric estimators.

The other approach is the semiparametric one, that is, we model $l$ parametrically. At the price of an additional model risk, parametric methods yield estimates that are always proper stable tail dependence functions. Moreover, they do not suffer from the curse of dimensionality. A large number of models have been proposed in the literature, allowing for various degrees of dependence and asymmetry, and new models continue to be invented; see [1, 20] for an overview of the most common ones.

In this paper, we propose an estimator based on the method of moments: given a parametric family $\{l_\theta : \theta \in \Theta\}$ with $\Theta \subseteq \mathbb{R}^p$ and a function $g : [0,1]^2 \to \mathbb{R}^p$, the moment estimator $\hat{\theta}_n$ is defined as the solution to the system of equations

$$\iint_{[0,1]^2} g(x,y) l_{\hat{\theta}_n}(x,y) \,\mathrm{d}x\,\mathrm{d}y = \iint_{[0,1]^2} g(x,y) \hat{l}_n(x,y) \,\mathrm{d}x\,\mathrm{d}y.$$

Here, $\hat{l}_n$ is the nonparametric estimator of $l$. Moreover, a comparison of the parametric and nonparametric estimators yields a goodness-of-fit test for the postulated model.

The method of moments estimator is to be contrasted with the maximum likelihood estimator in point process models for extremes [5, 17] or the censored likelihood approach proposed in [21, 23] and studied for single-parameter families in [14]. In parametric models, moment estimators yield consistent estimators, but often with a lower efficiency than the maximum likelihood estimator. However, as we shall see, the set of conditions required for the moment estimator is smaller, the conditions that remain to be imposed are much simpler and, most importantly, there are no restrictions whatsoever on the smoothness (or even on the existence) of the partial derivatives of $l$. Even for nonparametric estimators of $l$, asymptotic normality theorems require $l$ to be differentiable [6, 7, 15].

Such a degree of generality is needed if, for instance, the spectral measure underlying $l$ is discrete. In this case, there is no likelihood at all, so the maximum likelihood method breaks down. An example is the linear factor model $\boldsymbol{X} = \boldsymbol{\beta}\boldsymbol{F} + \boldsymbol{\varepsilon}$, where $\boldsymbol{X}$ and $\boldsymbol{\varepsilon}$ are $d \times 1$ random vectors, $\boldsymbol{F}$ is a $m \times 1$ random vector of factor variables and $\boldsymbol{\beta}$ is a constant $d \times m$ matrix of factor loadings. If the $m$ factor variables are mutually independent and if their common marginal tail is of Pareto type and heavier than those of the noise variables $\varepsilon_1, \ldots, \varepsilon_d$, then the spectral measure of the distribution of $\boldsymbol{X}$ is discrete with point masses determined by $\boldsymbol{\beta}$ and the tail index of the factor variables. The heuristic is that if $\boldsymbol{X}$ is far from the origin, then with high probability, it will be dominated by a single component of $\boldsymbol{F}$. Therefore, in the limit, there are only a finite number of directions for extreme outcomes of $\boldsymbol{X}$. Section 5 deals with a two-factor model of the above type, which gives rise to a discrete spectral measure concentrated on only two atoms. For more examples of factor models and further references, see [11].

The paper is organized as follows. Basic properties of stable tail dependence functions and spectral measures are reviewed in Section 2. The estimator and goodness-of-fit test statistic are defined in Section 3. Section 4 states the main results on the large-sample properties of the new procedures. In Section 5, the example of a spectral measure with two atoms is worked out and the finite-sample performance of the moment estimator is evaluated via simulations; Section 6 carries out the same program for the stable tail dependence functions of elliptical distributions. All proofs are deferred to Section 7.

## 2. Tail dependence

Let $(X, Y), (X_1, Y_1), \ldots, (X_n, Y_n)$ be independent random vectors in $\mathbb{R}^2$ with common continuous distribution function $F$ and marginal distribution functions $F_1$ and $F_2$. The central assumption in this paper is the existence, for all $(x, y) \in [0, \infty)^2$, of the limit $l$ in (1.1). Obviously, by the probability integral transform and the inclusion–exclusion formula, (1.1) is equivalent to the existence, for all $(x, y) \in [0, \infty)^2 \setminus \{(\infty, \infty)\}$, of the limit

$$\lim_{t \to 0} t^{-1} \mathbb{P}\{1 - F_1(X) \le tx, 1 - F_2(Y) \le ty\} = R(x, y), \tag{2.1}$$

so $R(x, \infty) = R(\infty, x) = x$. The functions $l$ and $R$ are related by $R(x, y) = x + y - l(x, y)$ for $(x, y) \in [0, \infty)^2$. Note that $R(1, 1)$ is the upper tail dependence coefficient.

If $C$ denotes the copula of $F$, that is, if $F(x, y) = C\{F_1(x), F_2(y)\}$, then (1.1) is equivalent to

$$\lim_{t \to 0} t^{-1} \{1 - C(1 - tx, 1 - ty)\} = l(x, y) \tag{2.2}$$

for all $x, y \ge 0$ and also to

$$\lim_{n \to \infty} C^n(u^{1/n}, v^{1/n}) = \exp\{-l(-\log u, -\log v)\} =: C_\infty(u, v)$$

for all $(u, v) \in (0, 1]^2$. The left-hand side in the previous display is the copula of the pair of componentwise maxima $(\max_{i=1,\ldots,n} X_i, \max_{i=1,\ldots,n} Y_i)$ and the right-hand side is the copula of a bivariate max-stable distribution. If, in addition, the marginal distribution functions $F_1$ and $F_2$ are in the max-domains of attraction of extreme value distributions $G_1$ and $G_2$, that is, if there exist normalizing sequences $a_n > 0, c_n > 0, b_n \in \mathbb{R}$ and $d_n \in \mathbb{R}$ such that $F_1^n(a_n x + b_n) \xrightarrow{d} G_1(x)$ and $F_2^n(c_n y + d_n) \xrightarrow{d} G_2(y)$, then actually

$$F^n(a_n x + b_n, c_n y + d_n) \xrightarrow{d} G(x,y) = C_\infty\{G_1(x), G_2(y)\},$$

that is, $F$ is in the max-domain of attraction of a bivariate extreme value distribution $G$ with marginals $G_1$ and $G_2$ and copula $C_\infty$. However, in this paper, we shall make no assumptions whatsoever on the marginal distributions $F_1$ and $F_2$, except for continuity.

Directly from the definition of $l$, it follows that $x \vee y \leq l(x,y) \leq x + y$ for all $(x,y) \in [0,\infty)^2$. Similarly, $0 \leq R(x,y) \leq x \wedge y$ for $(x,y) \in [0,\infty)^2$. Moreover, the functions $l$ and $R$ are homogeneous of order one: for all $(x,y) \in [0,\infty)^2$ and all $t \geq 0$,

$$l(tx, ty) = tl(x,y), \qquad R(tx, ty) = tR(x,y).$$

In addition, $l$ is convex and $R$ is concave. It can be shown that these requirements on $l$ (or, equivalently, $R$) are necessary and sufficient for $l$ to be a stable tail dependence function.

The following representation will be extremely useful: there exists a finite Borel measure $H$ on $[0,1]$, called *spectral* or *angular measure*, such that for all $(x,y) \in [0,\infty)^2$,

$$
\begin{aligned}
l(x,y) &= \int_{[0,1]} \max\{wx, (1-w)y\} H(\mathrm{d}w), \\
R(x,y) &= \int_{[0,1]} \min\{wx, (1-w)y\} H(\mathrm{d}w).
\end{aligned}
\tag{2.3}
$$

The identities $l(x,0) = l(0,x) = x$ for all $x \geq 0$ imply the following moment constraints for $H$:

$$\int_{[0,1]} w H(\mathrm{d}w) = \int_{[0,1]} (1-w) H(\mathrm{d}w) = 1. \tag{2.4}$$

Again, equation (2.4) constitutes a necessary and sufficient condition for $l$ in (2.3) to be a stable tail dependence function. For more details on multivariate extreme value theory, see, for instance, [1, 4, 8, 10, 13, 22].

## 3. Estimation and testing

Let $R_i^X$ and $R_i^Y$ be the rank of $X_i$ among $X_1, \ldots, X_n$ and the rank of $Y_i$ among $Y_1, \ldots, Y_n$, respectively, where $i = 1, \ldots, n$. Replacing $\mathbb{P}, F_1$ and $F_2$ on the left-hand side of (1.1) by their empirical counterparts, we obtain a nonparametric estimator for $l$. Estimators obtained in this way are

$$\hat{L}_n^1(x,y) := \frac{1}{k} \sum_{i=1}^n \mathbf{1}\{R_i^X > n + 1 - kx \text{ or } R_i^Y > n + 1 - ky\},$$

$$\hat{L}_n^2(x, y) := \frac{1}{k} \sum_{i=1}^n \mathbf{1}\{R_i^X \geq n + 1 - kx \text{ or } R_i^Y \geq n + 1 - ky\},$$

defined in [7] and [6, 15], respectively (here, $k \in \{1, \ldots, n\}$). The estimator we will use here is similar to those above and is defined by

$$\hat{l}_n(x, y) := \frac{1}{k} \sum_{i=1}^n \mathbf{1}\left\{R_i^X > n + \frac{1}{2} - kx \text{ or } R_i^Y > n + \frac{1}{2} - ky\right\}.$$

For finite samples, simulation experiments show that the latter estimator usually performs slightly better. The large-sample behaviors of the three estimators coincide, however, since $\hat{L}_n^1 \leq \hat{L}_n^2 \leq \hat{l}_n$ and, as $n \to \infty$,

$$\sup_{0 \leq x, y \leq 1} |\sqrt{k}(\hat{l}_n(x, y) - \hat{L}_n^1(x, y))| \leq \frac{2}{\sqrt{k}} \to 0, \tag{3.1}$$

where $k = k_n$ is an intermediate sequence, that is, $k \to \infty$ and $k/n \to 0$.

Assume that the stable tail dependence function $l$ belongs to some parametric family $\{l(\cdot, \cdot; \theta) : \theta \in \Theta\}$, where $\Theta \subset \mathbb{R}^p$, $p \geq 1$. (In the sequel, we will write $l(x, y; \theta)$ instead of $l_\theta(x, y)$.) Observe that this does not mean that $C$ (or $F$) belongs to a parametric family, that is, we have constructed a semiparametric model. Let $g : [0, 1]^2 \to \mathbb{R}^p$ be an integrable function such that $\varphi : \Theta \to \mathbb{R}^p$ defined by

$$\varphi(\theta) := \iint_{[0,1]^2} g(x, y) l(x, y; \theta) \, dx \, dy \tag{3.2}$$

is a homeomorphism between $\Theta^o$, the interior of the parameter space $\Theta$, and its image $\varphi(\Theta^o)$. For examples of the function $\varphi$, see Sections 5 and 6. Let $\theta_0$ denote the true parameter value and assume that $\theta_0 \in \Theta^o$.

The method of moments estimator $\hat{\theta}_n$ of $\theta_0$ is defined as the solution of

$$\iint_{[0,1]^2} g(x, y) \hat{l}_n(x, y) \, dx \, dy = \iint_{[0,1]^2} g(x, y) l(x, y; \hat{\theta}_n) \, dx \, dy = \varphi(\hat{\theta}_n),$$

that is,

$$\hat{\theta}_n := \varphi^{-1}\left(\iint_{[0,1]^2} g(x, y) \hat{l}_n(x, y) \, dx \, dy\right), \tag{3.3}$$

whenever the right-hand side is defined. For definiteness, if $\iint g\hat{l}_n \notin \varphi(\Theta^o)$, let $\hat{\theta}_n$ be some arbitrary, fixed value in $\Theta$.

Consider the goodness-of-fit testing problem, $\mathcal{H}_0 : l \in \{l(\cdot, \cdot; \theta) : \theta \in \Theta\}$ against $\mathcal{H}_a : l \notin \{l(\cdot, \cdot; \theta) : \theta \in \Theta\}$. We propose the test statistic

$$\iint_{[0,1]^2} \{\hat{l}_n(x, y) - l(x, y; \hat{\theta}_n)\}^2 \, dx \, dy, \tag{3.4}$$

with $\hat{\theta}_n$ as in (3.3). The null hypothesis is rejected for large values of the test statistic.

# 4. Results

The method of moments estimator is consistent for every intermediate sequence $k = k_n$ under minimal conditions on the model and the function $g$.

**Theorem 4.1 (Consistency).** *Let $g : [0, 1]^2 \to \mathbb{R}^p$ be integrable. If $\varphi$ in (3.2) is a homeomorphism between $\Theta^o$ and $\varphi(\Theta^o)$ and if $\theta_0 \in \Theta^o$, then as $n \to \infty$, $k \to \infty$ and $k/n \to 0$, the right-hand side of (3.3) is well defined with probability tending to 1 and $\hat{\theta}_n \overset{\mathbb{P}}{\to} \theta_0$.*

Denote by $W$ a mean-zero Wiener process on $[0, \infty]^2 \setminus \{(\infty, \infty)\}$ with covariance function

$$\mathbb{E}W(x_1, y_1)W(x_2, y_2) = R(x_1 \wedge x_2, y_1 \wedge y_2)$$

and for $x, y \in [0, \infty)$, let

$$W_1(x) := W(x, \infty), \qquad W_2(y) := W(\infty, y).$$

Further, for $(x, y) \in [0, \infty)^2$, let $R_1(x, y)$ and $R_2(x, y)$ be the right-hand partial derivatives of $R$ at the point $(x, y)$ with respect to the first and second coordinate, respectively. Since $R$ is concave, $R_1$ and $R_2$ defined in this way always exist, although they are discontinuous at points where $\frac{\partial}{\partial x} R(x, y)$ or $\frac{\partial}{\partial y} R(x, y)$ do not exist.

Finally, define the stochastic process $B$ on $[0, \infty)^2$ and the $p$-variate random vector $\tilde{B}$ by

$$B(x, y) = W(x, y) - R_1(x, y)W_1(x) - R_2(x, y)W_2(y),$$

$$\tilde{B} = \iint_{[0,1]^2} g(x, y)B(x, y)\, \mathrm{d}x\, \mathrm{d}y.$$

**Theorem 4.2 (Asymptotic normality).** *In addition to the conditions in Theorem 4.1, assume the following:*

(C1) *the function $\varphi$ is continuously differentiable in some neighborhood of $\theta_0$ and its derivative matrix $D_\varphi(\theta_0)$ is invertible;*

(C2) *there exists $\alpha > 0$ such that as $t \to 0$,*

$$t^{-1}\mathbb{P}\{1 - F_1(X) \le tx, 1 - F_2(Y) \le ty\} - R(x, y) = O(t^\alpha),$$

*uniformly on the set $\{(x, y) : x + y = 1, x \ge 0, y \ge 0\}$;*

(C3) *$k = k_n \to \infty$ and $k = o(n^{2\alpha/(1+2\alpha)})$ as $n \to \infty$.*

*Then*

$$\sqrt{k}(\hat{\theta}_n - \theta_0) \overset{d}{\to} D_\varphi(\theta_0)^{-1}\tilde{B}. \tag{4.1}$$

Note that condition (C2) is a second-order condition quantifying the speed of convergence in (2.1). Condition (C3) gives an upper bound on the speed with which $k$ can grow to infinity. This upper bound is related to the speed of convergence in (C2) and ensures that $\hat{\theta}_n$ is asymptotically unbiased.

The limiting distribution in (4.1) depends on the model and on the auxiliary function $g$. The optimal $g$ would be the one minimizing the asymptotic variance, but this minimization problem is typically difficult to solve. In the examples in Sections 5 and 6, the functions $g$ were chosen so as to simplify the calculations.

From the definition of the process $B$, it follows that the distribution of $\tilde{B}$ is $p$-variate normal with mean zero and covariance matrix

$$\Sigma(\theta_0) = \mathrm{Var}(\tilde{B}) = \iiiint_{[0,1]^4} g(x,y)g(u,v)^\top \sigma(x,y,u,v;\theta_0)\,\mathrm{d}x\,\mathrm{d}y\,\mathrm{d}u\,\mathrm{d}v, \tag{4.2}$$

where $\sigma$ is the covariance function of the process $B$, that is, for $\theta \in \Theta$,

$$\begin{aligned}
\sigma(x,y,u,v;\theta) &= \mathbb{E}B(x,y)B(u,v)\\
&= R(x \wedge u, y \wedge v;\theta) + R_1(x,y;\theta)R_1(u,v;\theta)(x \wedge u)\\
&\quad + R_2(x,y;\theta)R_2(u,v;\theta)(y \wedge v) - 2R_1(u,v;\theta)R(x \wedge u, y;\theta)\\
&\quad - 2R_2(u,v;\theta)R(x, y \wedge v;\theta) + 2R_1(x,y;\theta)R_2(u,v;\theta)R(x,v;\theta).
\end{aligned} \tag{4.3}$$

Denote by $H_\theta$ the spectral measure corresponding to $l(\cdot,\cdot;\theta)$. The following corollary allows the construction of confidence regions.

**Corollary 4.3.** *Under the assumptions of Theorem 4.2, if the map $\theta \mapsto H_\theta$ is weakly continuous at $\theta_0$ and if $\Sigma(\theta_0)$ is non-singular, then, as $n \to \infty$,*

$$k(\hat{\theta}_n - \theta_0)^\top D_\varphi(\hat{\theta}_n)^\top \Sigma(\hat{\theta}_n)^{-1} D_\varphi(\hat{\theta}_n)(\hat{\theta}_n - \theta_0) \xrightarrow{d} \chi_p^2.$$

Finally, we derive the limit distribution of the test statistic in (3.4).

**Theorem 4.4 (Test).** *Assume that the null hypothesis $\mathcal{H}_0$ holds and let $\theta_{\mathcal{H}_0}$ denote the true parameter. If*

(1) *for all $\theta_0 \in \Theta$ the conditions of Theorem 4.2 are satisfied (and hence $\Theta$ is open);*

(2) *on $\Theta$, the mapping $\theta \mapsto l(x,y;\theta)$ is differentiable for all $(x,y) \in [0,1]^2$ and its gradient is bounded in $(x,y) \in [0,1]^2$,*

*then*

$$\iint_{[0,1]^2} k(\hat{l}_n(x,y) - l(x,y;\hat{\theta}_n))^2\,\mathrm{d}x\,\mathrm{d}y$$

$$\xrightarrow{d} \iint_{[0,1]^2} (B(x,y) - D_{l(x,y;\theta)}(\theta_{\mathcal{H}_0})D_\varphi(\theta_{\mathcal{H}_0})^{-1}\tilde{B})^2\,\mathrm{d}x\,\mathrm{d}y$$

as $n \to \infty$, where $D_{l(x,y;\theta)}(\theta_{\mathcal{H}_0})$ is the gradient of $\theta \mapsto l(x,y;\theta)$ at $\theta_{\mathcal{H}_0}$.

# 5. Example 1: Two-point spectral measure

The *two-point spectral measure* is a spectral measure $H$ that is concentrated on only two points in $(0,1) \setminus \{1/2\}$ – call them $a$ and $1-b$. The moment conditions (2.4) imply that one of those points is less than $1/2$ and the other one is greater than $1/2$, and the masses on those points are determined by their locations. For definiteness, let $a \in (0, 1/2)$ and $1 - b \in (1/2, 1)$, so the parameter vector $\theta = (a,b)$ takes values in the square $\Theta = (0, 1/2)^2$. The masses assigned to $a$ and $1-b$ are

$$q := H(\{a\}) = \frac{1 - 2b}{1 - a - b} \quad \text{and} \quad 2 - q = H(\{1-b\}) = \frac{1 - 2a}{1 - a - b}.$$

This model is also known as the *natural model* and was first described by Tiago de Oliveira [24, 25].

By (2.3), the corresponding stable tail dependence function is

$$l(x, y; a, b) = q \max\{ax, (1-a)y\} + (2 - q) \max\{(1-b)x, by\}.$$

The partial derivatives of $l$ with respect to $x$ and $y$ are

$$\frac{\partial l(x,y;a,b)}{\partial x} = \begin{cases} 1, & \text{if } y < \dfrac{a}{1-a}x, \\ (1-b)(2-q), & \text{if } \dfrac{a}{1-a}x < y < \dfrac{1-b}{b}x, \\ 0, & \text{if } y > \dfrac{1-b}{b}x \end{cases}$$

and $(\partial/\partial y)l(x,y;a,b) = (\partial/\partial y)l(y,x;b,a)$. Note that the partial derivatives do not exist on the lines $y = \frac{a}{1-a}x$ and $y = \frac{1-b}{b}x$. The same is true for the partial derivatives of $R$. As a consequence, the maximum likelihood method is not applicable and the asymptotic normality of the nonparametric estimator breaks down. However, the method of moments estimator can still be used since, in Theorem 4.2, no smoothness assumptions whatsoever are made on $l$.

As explained in the Introduction, discrete spectral measures arise whenever extremes are determined by a finite number of independent, heavy-tailed factors. Specifically, let the random vector $(X, Y)$ be given by

$$(X, Y) = (\alpha Z_1 + (1 - \alpha)Z_2 + \varepsilon_1, (1 - \beta)Z_1 + \beta Z_2 + \varepsilon_2), \tag{5.1}$$

where $0 < \alpha < 1$ and $0 < \beta < 1$ are coefficients and where $Z_1$, $Z_2$, $\varepsilon_1$ and $\varepsilon_2$ are independent random variables satisfying the following conditions: there exist $\nu > 0$ and a slowly varying function $L$ such that $\mathbb{P}(Z_i > z) = z^{-\nu}L(z)$ for some $\nu > 0$, $i = 1, 2$; $\mathbb{P}(\varepsilon_j > z)/\mathbb{P}(Z_1 > z) \to 0$ as $z \to \infty$, $j = 1, 2$. (Recall that a positive, measurable function $L$ defined in a neighborhood of infinity is called *slowly varying* if $L(yz)/L(z) \to 1$ as $z \to \infty$ for all $y > 0$.) Straightforward, but lengthy, computations show that the spectral measure of the random vector $(X, Y)$ is a two-point spectral measure having masses $q$ and $2 - q$ at the points $a$ and $1 - b$, where

$$q := \frac{(1-\alpha)^\nu}{\alpha^\nu + (1-\alpha)^\nu} + \frac{\beta^\nu}{\beta^\nu + (1-\beta)^\nu},$$

$$a := \frac{(1-\alpha)^\nu}{\alpha^\nu + (1-\alpha)^\nu} q^{-1}, \qquad 1 - b := \frac{\alpha^\nu}{\alpha^\nu + (1-\alpha)^\nu} (2-q)^{-1}.$$

Write $\Delta = \{(x,y) \in [0,1]^2 : x + y \le 1\}$ and let $\mathbf{1}_\Delta$ be its indicator function. The function $g_\Delta : [0,1]^2 \to \mathbb{R}^2$ defined by $g_\Delta(x,y) = \mathbf{1}_\Delta(x,y)(x,y)^\top$ is obviously integrable and the function $\varphi$ in (3.2) is given by

$$\varphi(a,b) = \iint_\Delta (x,y)^\top l(x,y;a,b) \, \mathrm{d}x \, \mathrm{d}y = (J(a,b), K(b,a))^\top,$$

where $K(a,b) = J(b,a)$ and

$$J(a,b) = \tfrac{1}{24} \{(2ab - a - b)(b - a + 1) + a(b-1) + 3\}.$$

Nonparametric estimators of $J$ and $K$ are given by

$$(\hat{J}_n, \hat{K}_n) = \iint_\Delta (x,y)^\top \hat{l}_n(x,y) \, \mathrm{d}x \, \mathrm{d}y$$

and the method of moment estimators $(\hat{a}_n, \hat{b}_n)$ are defined as the solutions to the equations

$$(\hat{J}_n, \hat{K}_n) = (J(\hat{a}_n, \hat{b}_n), K(\hat{a}_n, \hat{b}_n)).$$

Due to the explicit nature of the functions $J$ and $K$, these equations can be simplified: if we denote $c_{J,n} := 3(8\hat{J}_n - 1)$ and $c_{K,n} := 3(8\hat{K}_n - 1)$, the estimator $\hat{b}_n$ of $b$ will be a solution of the quadratic equation

$$3(2c_{J,n} + 2c_{K,n} + 3)b^2 + 3(-5c_{J,n} + c_{K,n} - 3)b + 3c_{J,n} - 6c_{K,n} - (c_{J,n} + c_{K,n})^2 = 0$$

that falls into the interval $(0, 1/2)$ and the estimator of $a$ is

$$\hat{a}_n = \frac{3\hat{b}_n + c_{J,n} + c_{K,n}}{6\hat{b}_n - 3}.$$

In the simulations, we used the following models:

(i) $Z_1, Z_2 \sim \text{Fréchet}(1)$, so $\nu = 1$, and $\varepsilon_1, \varepsilon_2 \sim N(0,1)$ (Figures 1, 2, 3);
(ii) $Z_1, Z_2 \sim t_2$, so $\nu = 1/2$, and $\varepsilon_1, \varepsilon_2 \sim N(0, 0.5^2)$ (Figures 4, 5, 6).

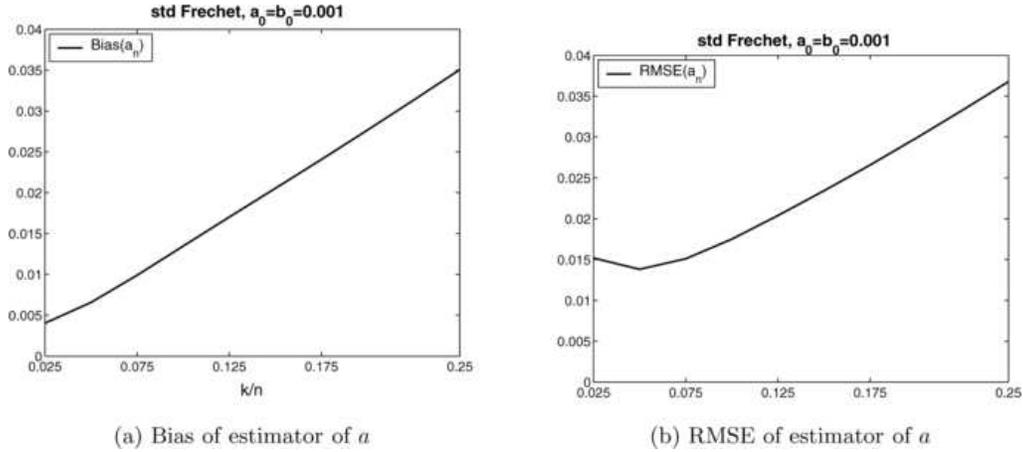

(a) Bias of estimator of $a$                    (b) RMSE of estimator of $a$

**Figure 1.** Model (5.1) with $Z_1, Z_2 \sim$ Fréchet(1), $\varepsilon_1, \varepsilon_2 \sim N(0,1)$, $a_0 = b_0 = 0.001$.

The figures show the bias and the root mean squared error (RMSE) of $\hat{a}_n$ and $\hat{b}_n$ for 1000 samples of size $n = 1000$. The method of moments estimator performs well in general. We see a very good behavior when $a_0 = b_0 \approx 0$. Of course, the heavier the tail of $Z_1$, the better the performance of the estimator.

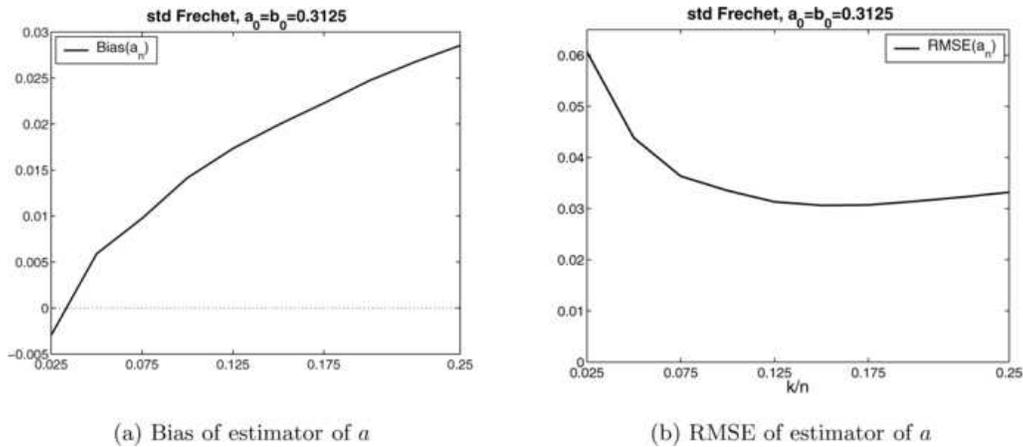

(a) Bias of estimator of $a$                    (b) RMSE of estimator of $a$

**Figure 2.** Model (5.1) with $Z_1, Z_2 \sim$ Fréchet(1), $\varepsilon_1, \varepsilon_2 \sim N(0,1)$, $a_0 = b_0 = 0.3125$.

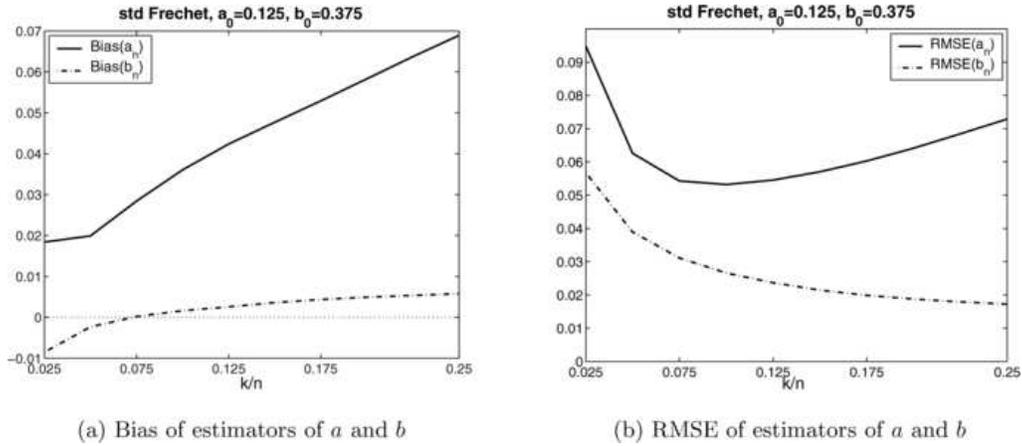

(a) Bias of estimators of $a$ and $b$      (b) RMSE of estimators of $a$ and $b$

**Figure 3.** Model (5.1) with $Z_1, Z_2 \sim$ Fréchet(1), $\varepsilon_1, \varepsilon_2 \sim N(0,1)$, $a_0 = 0.125$, $b_0 = 0.375$.

## 6. Example 2: Parallel meta-elliptical model

A random vector $(X, Y)$ is said to be *elliptically distributed* if it satisfies the distributional equality

$$(X, Y)^\top \stackrel{d}{=} \boldsymbol{\mu} + Z\boldsymbol{A}\boldsymbol{U}, \tag{6.1}$$

where $\boldsymbol{\mu}$ is a $2 \times 1$ column vector, $Z$ is a positive random variable called the generating random variable, $\boldsymbol{A}$ is a $2 \times 2$ matrix such that $\boldsymbol{\Sigma} = \boldsymbol{A}\boldsymbol{A}^\top$ is of full rank and $\boldsymbol{U}$ is a two-dimensional random vector independent of $Z$ and uniformly distributed on the unit circle $\{(x, y) \in \mathbb{R}^2 : x^2 + y^2 = 1\}$. Under the above assumptions, the matrix $\boldsymbol{\Sigma}$ can be written as

$$\boldsymbol{\Sigma} = \begin{pmatrix} \sigma^2 & \rho\sigma v \\ \rho\sigma v & v^2 \end{pmatrix}, \tag{6.2}$$

where $\sigma > 0$, $v > 0$ and $-1 < \rho < 1$. The special case $\rho = 0$ yields the subclass of *parallel elliptical distributions*.

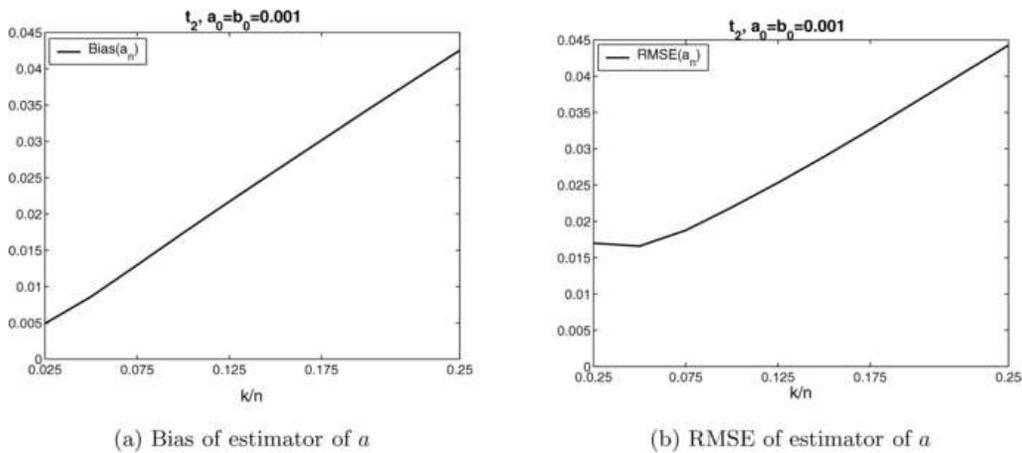

(a) Bias of estimator of $a$      (b) RMSE of estimator of $a$

**Figure 4.** Model (5.1) with $Z_1, Z_2 \sim t_2$, $\varepsilon_1, \varepsilon_2 \sim N(0, 0.5^2)$, $a_0 = b_0 = 0.001$.

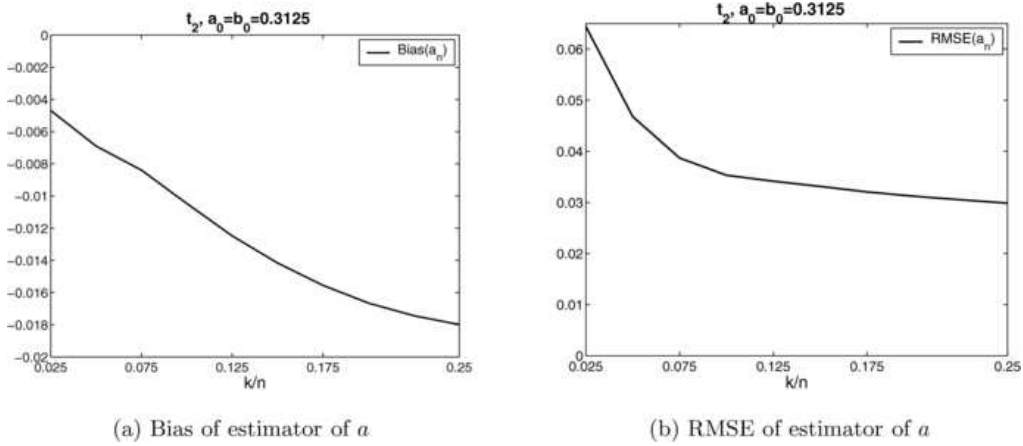

(a) Bias of estimator of $a$                (b) RMSE of estimator of $a$

**Figure 5.** Model (5.1) with $Z_1, Z_2 \sim t_2$, $\varepsilon_1, \varepsilon_2 \sim N(0, 0.5^2)$, $a_0 = b_0 = 0.3125$.

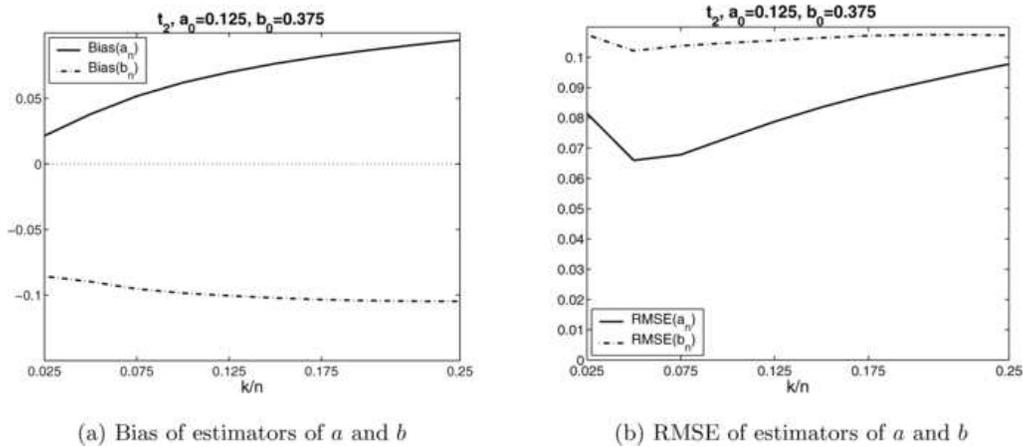

(a) Bias of estimators of $a$ and $b$                (b) RMSE of estimators of $a$ and $b$

**Figure 6.** Model (5.1) with $Z_1, Z_2 \sim t_2$, $\varepsilon_1, \varepsilon_2 \sim N(0, 0.5^2)$, $a_0 = 0.125$, $b_0 = 0.375$.

By [16], the distribution of $Z$ satisfies $\mathbb{P}(Z > z) = z^{-\nu} L(z)$ with $\nu > 0$ and $L$ slowly varying if and only if the distribution of $(X, Y)$ is (multivariate) regularly varying with the same index. Under this assumption, the function $R$ of the distribution of $(X, Y)$ was derived in [18]. In case $\rho = 0$, the formula specializes to

$$R(x, y; \nu) = \frac{x \int_{f(x,y;\nu)}^{\pi/2} (\cos \phi)^{\nu} \, \mathrm{d}\phi + y \int_0^{f(x,y;\nu)} (\sin \phi)^{\nu} \, \mathrm{d}\phi}{\int_{-\pi/2}^{\pi/2} (\cos \phi)^{\nu} \, \mathrm{d}\phi} \tag{6.3}$$

with $f(x, y; \nu) = \arctan\{(x/y)^{1/\nu}\}$. Hence, the class of stable tail dependence functions belonging to parallel elliptical vectors with regularly varying generating random variables forms a one-dimensional parametric family indexed by the index of regular variation $\nu \in (0, \infty) = \Theta$ of $Z$. We will call the corresponding stable tail dependence functions $l$ *parallel elliptical*.

In [9], *meta-elliptical distributions* are defined as the distributions of random vectors of the form $(s(X), t(Y))$, where the distribution of $(X, Y)$ is elliptical and $s$ and $t$ are increasing functions. In other words, a distribution is meta-elliptical if and only if its copula is that of an elliptical distribution. Such copulas are called meta-elliptical in [12] (note that a copula, as a distribution function on the unit square, cannot be elliptical in the sense of (6.1)). Since a stable tail dependence function $l$ of a bivariate distribution $F$ is only determined by $F$ through its copula $C$ (see (2.2)), the results in the preceding paragraph continue to hold for meta-elliptical distributions. In the case $\rho = 0$, we will speak of parallel meta-elliptical distributions. In the case where the generating random variable $Z$ is regularly varying with index $\nu$, the function $R$ is given by (6.3).

For parallel meta-elliptical distributions, the second-order condition (C2) in Theorem 4.2 can be checked via second-order regular variation of $Z$.

**Lemma 6.1.** *Let $F$ be a parallel meta-elliptical distribution with generating random variable $Z$. If there exist $\nu > 0$, $\beta < 0$ and a function $A(t) \to 0$ of constant sign near infinity such that*

$$\lim_{t \to \infty} \frac{\mathbb{P}(Z > tx)/\mathbb{P}(Z > t) - x^{-\nu}}{A(t)} = x^{-\nu} \frac{x^{\beta} - 1}{\beta}, \tag{6.4}$$

*then condition* (C2) *in Theorem 4.2 holds for every $\alpha \in (0, -\beta/\nu)$.*

Note that although the generating random variable is only defined up to a multiplicative constant, condition (6.4) does makes sense: that is, if (6.4) holds for a random variable $Z$, then it also holds for $cZ$ with $c > 0$, for the same constants $\nu$ and $\beta$ and for the rate function $A^*(t) := A(t/c)$. Note that $|A|$ is necessarily regularly varying with index $\beta$; see [2], equation (3.0.3).

Now, assume that $(X_1, Y_1), \ldots, (X_n, Y_n)$ is a random sample from a bivariate distribution $F$ with parallel elliptical stable tail dependence function $l$, that is, $l \in \{l(\cdot, \cdot; \nu) : \nu \in (0, \infty)\}$, where $l(x, y; \nu) = x + y - R(x, y; \nu)$ and $R(x, y; \nu)$ is as in (6.3). We will apply the method of moments to estimate the parameter $\nu$. Since $l$ is defined by a limit relation, our

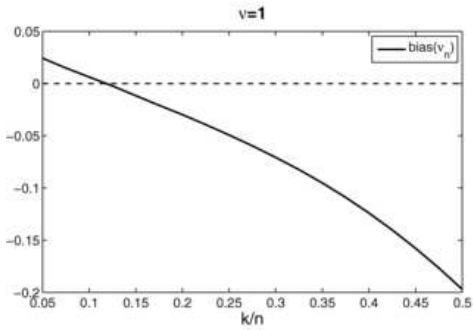

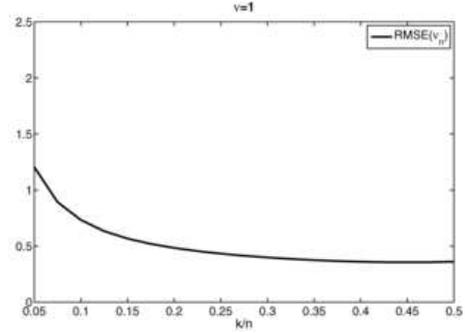

(a) Bias of estimator of $\nu$

(b) RMSE of estimator of $\nu$

**Figure 7.** Estimation of $\nu = 1$ in the bivariate Cauchy model.

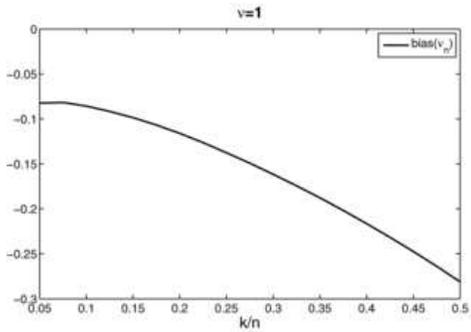

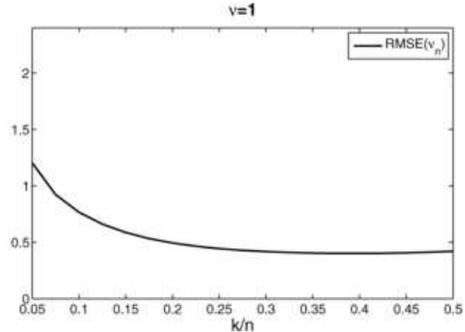

(a) Bias of estimator of $\nu$

(b) RMSE of estimator of $\nu$

**Figure 8.** Estimation of $\nu = 1$ in the model $(X_1, Y_1)^\top = Z\boldsymbol{U}$, where $Z$ is Fréchet(1).

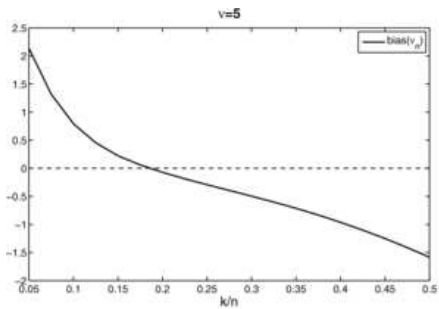

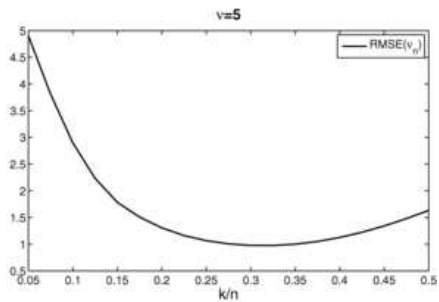

(a) Bias of estimator of $\nu$

(b) RMSE of estimator of $\nu$

**Figure 9.** Estimation of $\nu = 5$ in the model $(X_1, Y_1)^\top = Z\boldsymbol{U}$, where $Z$ is Fréchet(5).

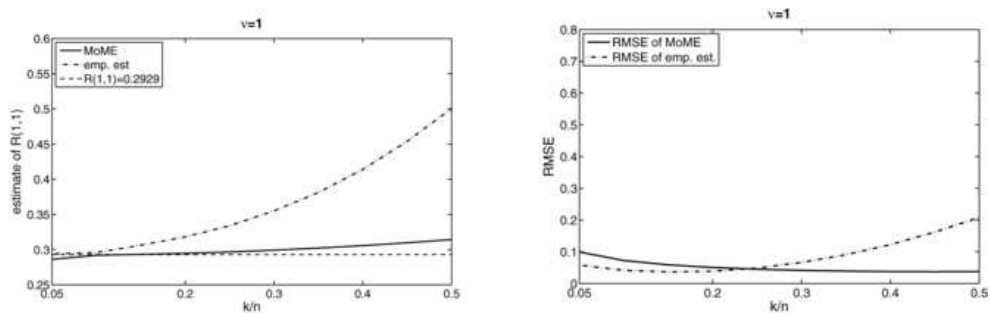

(a) the mean over 1000 replications of the estimate of $R(1,1;1)$

(b) the RMSE of the estimator of $R(1,1;1)$

**Figure 10.** Estimation of $R(1,1;1)$ in the bivariate Cauchy model.

assumption on $F$ is weaker than the assumption that $F$ is parallel meta-elliptical with regularly varying $Z$, which, as explained above, is, in turn, weaker than the assumption that $F$ itself is parallel elliptical with regularly varying $Z$. The problem of estimating the $R$ for elliptical distributions was addressed in [18] and for meta-elliptical distributions was addressed in [19].

We simulated 1000 random samples of size $n = 1000$ from models for which the assumptions of Theorem 4.2 hold and which have the function $R(\cdot, \cdot; \nu)$ as in (6.3), with $\nu \in \{1, 5\}$. The three models we used are of the type $(X_1, Y_1)^\top = Z\boldsymbol{U}$. In the first model, the generating random variable $Z$ is such that $\mathbb{P}(Z > z) = (1 + z^2)^{-1/2}$ for $z \geq 0$, that is, the first model is the bivariate Cauchy ($\nu = 1$). In the other two models, $Z$ is Fréchet($\nu$) with $\nu \in \{1, 5\}$.

Figures 7 to 9 show the bias and the RMSE of the moment estimator of $\nu$. The auxiliary function $g : [0, 1]^2 \to \mathbb{R}$ is $g(x, y) = \mathbf{1}(x + y \leq 1)$. For comparison, Figures 10 and 11 show the plots of the means and RMSE of the parametric and nonparametric estimates $R(1, 1; \hat{\nu}_n)$ and $\hat{R}_n(1, 1) = 2 - \hat{l}_n(1, 1)$ of the upper tail dependence coefficient $R(1, 1)$. We can see that the method of moments estimator of the upper tail dependence coefficient $R(1, 1; \nu)$ performs well. In particular, it is much less sensitive to the choice of $k$ than the nonparametric estimator.

## 7. Proofs

**Proof of Theorem 4.1.** First, note that

$$\left| \iint_{[0,1]^2} g(x, y) \hat{l}_n(x, y) \, \mathrm{d}x \, \mathrm{d}y - \iint_{[0,1]^2} g(x, y) l(x, y; \theta_0) \, \mathrm{d}x \, \mathrm{d}y \right|$$

$$\leq \sup_{0 \leq x, y \leq 1} |\hat{l}_n(x,y) - l(x,y;\theta_0)| \iint_{[0,1]^2} |g(x,y)| \, dx \, dy.$$

The second term is finite by assumption and

$$\sup_{0 \leq x, y \leq 1} |\hat{l}_n(x,y) - l(x,y;\theta_0)| \xrightarrow{\mathbb{P}} 0$$

by (3.1) and [15], Theorem 1; see also [6]. Therefore, as $n \to \infty$,

$$\iint_{[0,1]^2} g(x,y)\hat{l}_n(x,y) \, dx \, dy \xrightarrow{\mathbb{P}} \iint_{[0,1]^2} g(x,y)l(x,y;\theta_0) \, dx \, dy = \varphi(\theta_0).$$

Since $\varphi(\theta_0) \in \varphi(\Theta^o)$, which is open, and since $\varphi^{-1}$ is continuous at $\varphi(\theta_0)$ by assumption, we can apply the function $\varphi^{-1}$ on both sides of the previous limit relation so that, by the continuous mapping theorem, we indeed have $\hat{\theta}_n \xrightarrow{\mathbb{P}} \theta_0$. $\qquad \square$

For the proof of Theorem 4.2, we will need the following lemma, the proof of which follows from [8], Lemma 6.2.1.

**Lemma 7.1.** *The function R in (2.3) is differentiable at $(x,y) \in (0,\infty)^2$ if $H(\{z\}) = 0$ with $z = y/(x+y)$. In that case, the gradient of R is given by $(R_1(x,y), R_2(x,y))^\top$, where*

$$R_1(x,y) = \int_0^z wH(dw), \qquad R_2(x,y) = \int_z^1 (1-w)H(dw). \tag{7.1}$$

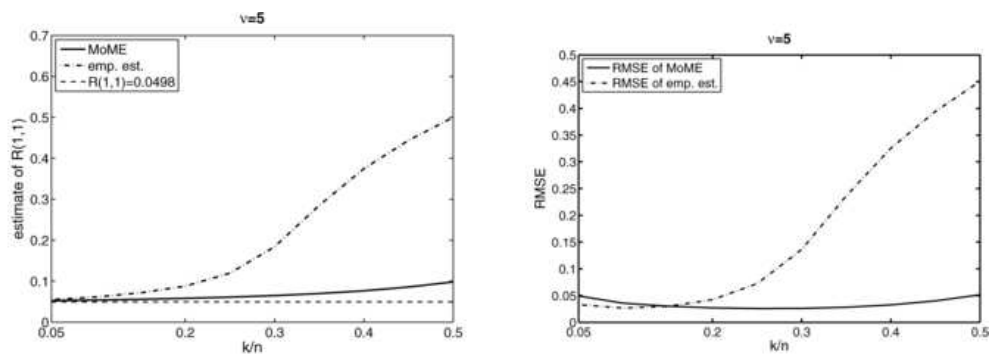

(a) the mean over 1000 replications of the estimate of $R(1,1;5)$

(b) the RMSE of the estimator of $R(1,1;5)$

**Figure 11.** Estimation of $R(1,1;5)$ in the model $(X_1, Y_1)^\top = Z\boldsymbol{U}$, where $Z$ is Fréchet(5).

For $i = 1, \ldots, n$, let $U_i := 1 - F_1(X_i)$ and $V_i := 1 - F_2(Y_i)$. Let $Q_{1n}$ and $Q_{2n}$ denote the empirical quantile functions of $(U_1, \ldots, U_n)$ and $(V_1, \ldots, V_n)$, respectively, that is,

$$Q_{1n}\left(\frac{kx}{n}\right) = U_{\lceil kx \rceil : n}, \qquad Q_{2n}\left(\frac{ky}{n}\right) = V_{\lceil ky \rceil : n},$$

where $U_{1:n} \leq \cdots \leq U_{n:n}$ and $V_{1:n} \leq \cdots \leq V_{n:n}$ are the order statistics and where $\lceil a \rceil$ is the smallest integer not smaller than $a$. Define

$$S_{1n}(x) := \frac{n}{k} Q_{1n}\left(\frac{kx}{n}\right), \qquad S_{2n}(y) := \frac{n}{k} Q_{2n}\left(\frac{ky}{n}\right)$$

and

$$\hat{R}_n^1(x,y) := \frac{1}{k} \sum_{i=1}^n \mathbf{1}\left\{ U_i < \frac{k}{n} S_{1n}(x), V_i < \frac{k}{n} S_{2n}(y) \right\},$$

$$= \frac{1}{k} \sum_{i=1}^n \mathbf{1}\{ U_i < U_{\lceil kx \rceil : n}, V_i < V_{\lceil ky \rceil : n} \},$$

$$= \frac{1}{k} \sum_{i=1}^n \mathbf{1}\{ R_i^X > n + 1 - kx, R_i^Y > n + 1 - ky \},$$

$$R_n(x,y) := \frac{n}{k} \mathbb{P}\left( U_1 \leq \frac{kx}{n}, V_1 \leq \frac{ky}{n} \right),$$

$$T_n(x,y) := \frac{1}{k} \sum_{i=1}^n \mathbf{1}\left\{ U_i < \frac{kx}{n}, V_i < \frac{ky}{n} \right\}.$$

Further, note that

$$\hat{R}_n^1(x,y) = T_n(S_{1n}(x), S_{2n}(y)).$$

Write $v_n(x,y) = \sqrt{k}(T_n(x,y) - R_n(x,y))$, $v_{n,1}(x) := v_n(x,\infty)$ and $v_{n,2}(y) := v_n(\infty, y)$. From [7], Proposition 3.1 we get

$$(v_n(x,y), x,y \in [0,1]; v_{n,1}(x), x \in [0,1]; v_{n,2}(y), y \in [0,1])$$

$$\xrightarrow{d} (W(x,y), x,y \in [0,1]; W_1(x), x \in [0,1]; W_2(y), y \in [0,1]),$$

in the topology of uniform convergence, as $n \to \infty$. Invoking the Skorokhod construction (see, e.g., [27]) we get a new probability space containing all $\tilde{v}_n, \tilde{v}_{n,1}, \tilde{v}_{n,2}, \tilde{W}, \tilde{W}_1, \tilde{W}_2$ for which it holds that

$$(\tilde{v}_n, \tilde{v}_{n,1}, \tilde{v}_{n,2}) \stackrel{d}{=} (v_n, v_{n,1}, v_{n,2}),$$

$$(\tilde{W}, \tilde{W}_1, \tilde{W}_2) \stackrel{d}{=} (W, W_1, W_2)$$

as well as

$$\sup_{0 \le x, y \le 1} |\tilde{v}_n(x,y) - \tilde{W}(x,y)| \overset{a.s.}{\to} 0,$$

$$\sup_{0 \le x \le 1} |\tilde{v}_{n,j}(x) - \tilde{W}_j(x)| \overset{a.s.}{\to} 0, \qquad j = 1, 2.$$

We will work on this space from now on, but keep the old notation (without tildes). The following consequence of the above and Vervaat's lemma [28] will be useful

$$\sup_{0 \le x \le 1} |\sqrt{k}(S_{jn}(x) - x) + W_j(x)| \overset{a.s.}{\to} 0, \qquad j = 1, 2. \tag{7.2}$$

**Proof of Theorem 4.2.** In this proof, we will write $l(x,y)$ and $R(x,y)$ instead of $l(x,y;\theta_0)$ and $R(x,y;\theta_0)$, respectively.

First, we will show that as $n \to \infty$,

$$\left| \sqrt{k} \left( \iint_{[0,1]^2} g(x,y) \hat{L}_n^1(x,y) \, dx \, dy - \varphi(\theta_0) \right) + \tilde{B} \right| \overset{\mathbb{P}}{\to} 0. \tag{7.3}$$

Since, for each $x, y \in (0,1]$,

$$(\hat{L}_n^1 + \hat{R}_n^1)(x,y) = \frac{\lceil kx \rceil + \lceil ky \rceil - 2}{k}$$

almost surely, from

$$\left| \frac{\lceil kx \rceil + \lceil ky \rceil - 2}{k} - x - y \right| \le \frac{2}{k},$$

it follows that

$$\left| \sqrt{k} \left( \iint_{[0,1]^2} g(x,y) \hat{L}_n^1(x,y) \, dx \, dy - \iint_{[0,1]^2} g(x,y) l(x,y) \, dx \, dy \right) \right.$$

$$\left. + \sqrt{k} \left( \iint_{[0,1]^2} g(x,y) \hat{R}_n^1(x,y) \, dx \, dy - \iint_{[0,1]^2} g(x,y) R(x,y) \, dx \, dy \right) \right|$$

$$= \left| \iint_{[0,1]^2} g(x,y) \sqrt{k} \left( \frac{\lceil kx \rceil + \lceil ky \rceil - 2}{k} - x - y \right) dx \, dy \right| = O\left( \frac{1}{\sqrt{k}} \right)$$

almost surely. Hence, to show (7.3), we will prove

$$\left| \iint_{[0,1]^2} g(x,y) \sqrt{k} (\hat{R}_n^1(x,y) - R(x,y)) \, dx \, dy - \tilde{B} \right| \overset{\mathbb{P}}{\to} 0. \tag{7.4}$$

First, we write

$$\sqrt{k} (\hat{R}_n^1(x,y) - R(x,y)) = \sqrt{k} (\hat{R}_n^1(x,y) - R_n(S_{1n}(x), S_{2n}(y)))$$

$$+ \sqrt{k}(R_n(S_{1n}(x), S_{2n}(y)) - R(S_{1n}(x), S_{2n}(y)))$$
$$+ \sqrt{k}(R(S_{1n}(x), S_{2n}(y)) - R(x, y)).$$

From the assumption on integrability of $g$ and the proof of [7], Theorem 2.2, page 2003, we get

$$\iint_{[0,1]^2} |g(x,y)| |\sqrt{k}(\hat{R}_n^1(x,y) - R_n(S_{1n}(x), S_{2n}(y))) - W(x,y)| \, \mathrm{d}x \, \mathrm{d}y$$

$$\leq \sup_{0 \leq x, y \leq 1} |\sqrt{k}(\hat{R}_n^1(x,y) - R_n(S_{1n}(x), S_{2n}(y))) - W(x,y)| \iint_{[0,1]^2} |g(x,y)| \, \mathrm{d}x \, \mathrm{d}y \quad (7.5)$$

$$\xrightarrow{\mathbb{P}} 0$$

and, by conditions (C2) and (C3),

$$\iint_{[0,1]^2} |g(x,y)| |\sqrt{k}(R_n(S_{1n}(x), S_{2n}(y)) - R(S_{1n}(x), S_{2n}(y)))| \, \mathrm{d}x \, \mathrm{d}y$$

$$\leq \sup_{0 \leq x, y \leq 1} |\sqrt{k}(R_n(S_{1n}(x), S_{2n}(y)) - R(S_{1n}(x), S_{2n}(y)))| \iint_{[0,1]^2} |g(x,y)| \, \mathrm{d}x \, \mathrm{d}y \quad (7.6)$$

$$\xrightarrow{\mathbb{P}} 0.$$

Take $\omega$ in the Skorokhod probability space introduced above such that $\sup_{0 \leq x \leq 1} |W_1(x)|$ and $\sup_{0 \leq y \leq 1} |W_2(y)|$ are finite and (7.2) holds. For such $\omega$, we will show, by means of dominated convergence, that

$$\iint_{[0,1]^2} |g(x,y)| |\sqrt{k}(R(S_{1n}(x), S_{2n}(y)) - R(x,y))$$
$$+ R_1(x,y)W_1(x) + R_2(x,y)W_2(y)| \, \mathrm{d}x \, \mathrm{d}y \to 0. \quad (7.7)$$

(i) *Pointwise convergence of the integrand to zero for almost all* $(x,y) \in [0,1]^2$. Convergence in $(x,y)$ follows from (7.2), provided $R(x,y)$ is differentiable. The set of points in which this might fail is, by Lemma 7.1, equal to

$$D_R := \left\{ (x,y) \in [0,1]^2 : H(\{z\}) > 0, z = \frac{y}{x+y} \right\}.$$

Since $H$ is a finite measure, there can be at most countably many $z$ for which $H(\{z\}) > 0$. The set $D_R$ is then a union of at most countably many lines through the origin and hence has Lebesgue measure zero.

(ii) *The domination of the integrand for all* $(x,y) \in [0,1]^2$. Comparing (7.1) and the moment conditions (2.4), we see that for all $(x,y) \in [0,1]^2$, it holds that $|R_1(x,y)| \leq 1$ and $|R_2(x, y)| \leq 1$. Hence, for all $(x, y) \in [0, 1]^2$,

$$|g(x, y)||\sqrt{k}(R(R(S_{1n}(x), S_{2n}(y)) - R(x, y)) + R_1(x, y)W_1(x) + R_2(x, y)W_2(y)|$$

$$\leq |g(x, y)|(\sqrt{k}|R(S_{1n}(x), S_{2n}(y)) - R(x, y)| + |W_1(x)| + |W_2(y)|).$$

We will show that the right-hand side in the above inequality is less than or equal to $M|g(x, y)|$ for all $(x, y) \in [0, 1]^2$ and some positive constant $M$ (depending on $\omega$). For that purpose, we prove that

$$\sup_{0 \leq x, y \leq 1} \sqrt{k}|R(S_{1n}(x), S_{2n}(y)) - R(x, y)| = O(1).$$

The representation (2.1) implies that for all $x, x_1, x_2, y, y_1, y_2 \in [0, 1]$,

$$|R(x_1, y) - R(x_2, y)| \leq |x_1 - x_2|,$$

$$|R(x, y_1) - R(x, y_2)| \leq |y_1 - y_2|.$$

By these inequalities and (7.2), we now have

$$\sup_{0 \leq x, y \leq 1} \sqrt{k}|R(S_{1n}(x), S_{2n}(y)) - R(x, y)|$$

$$\leq \sup_{0 \leq x, y \leq 1} \sqrt{k}|R(S_{1n}(x), S_{2n}(y)) - R(S_{1n}(x), y)| + \sup_{0 \leq x, y \leq 1} \sqrt{k}|R(S_{1n}(x), y) - R(x, y)|$$

$$\leq \sup_{0 \leq x \leq 1} \sqrt{k}|S_{1n}(x) - x| + \sup_{0 \leq y \leq 1} \sqrt{k}|S_{2n}(y) - y|$$

$$= O(1).$$

Recalling that $\sup_{0 \leq x \leq 1} |W_1(x)|$ and $\sup_{0 \leq y \leq 1} |W_2(y)|$ are finite completes the proof of domination and hence the proof of (7.7).

Combining (7.5), (7.6) and (7.7), we get (7.4) and therefore also (7.3). Property (3.1) provides us with a statement analogous to (7.3), but with $\hat{L}_n^1$ replaced by $\hat{l}_n$. That is, we have

$$\left| \sqrt{k} \left( \iint_{[0,1]^2} g(x, y) \hat{l}_n(x, y) \, dx \, dy - \varphi(\theta_0) \right) + \tilde{B} \right| \xrightarrow{\mathbb{P}} 0. \tag{7.8}$$

Using condition (C1) and the inverse mapping theorem, we get that $\varphi^{-1}$ is continuously differentiable in a neighborhood of $\varphi(\theta_0)$ and $D_{\varphi^{-1}}(\varphi(\theta_0))$ is equal to $D_\varphi(\theta_0)^{-1}$. By a routine argument, using the delta method (see, e.g., Theorem 3.1 in [26]), (7.8) implies that

$$\sqrt{k}(\hat{\theta}_n - \theta_0) \xrightarrow{\mathbb{P}} -D_\varphi(\theta_0)^{-1}\tilde{B}$$

and since $\tilde{B}$ is mean-zero normally distributed ($\tilde{B} \overset{d}{=} -\tilde{B}$),

$$\sqrt{k}(\hat{\theta}_n - \theta_0) \xrightarrow{d} D_\varphi(\theta_0)^{-1}\tilde{B}. \qquad \square$$

**Lemma 7.2.** *Let $H_\theta$ be the spectral measure and $\Sigma(\theta)$ the covariance matrix in (4.2). If the mapping $\theta \mapsto H_\theta$ is weakly continuous at $\theta_0$, then $\theta \mapsto \Sigma(\theta)$ is continuous at $\theta_0$.*

**Proof.** Let $\theta_n \to \theta_0$. In view of the expression for $\Sigma(\theta)$ in (4.2) and (4.3), the assumption that $g$ is integrable and the fact that $R$, $|R_1|$ and $|R_2|$ are bounded by 1 for all $\theta$ and $(x, y) \in [0, 1]^2$, it suffices to show that $R(x, y; \theta_n) \to R(x, y; \theta)$ and $R_i(x, y; \theta_n) \to R_i(x, y; \theta)$ for $i = 1, 2$ and for almost all $(x, y) \in [0, 1]^2$.

Convergence of $R$ for all $(x, y) \in [0, 1]^2$ follows directly from the representation of $R$ in terms of $H$ in (2.3) and the definition of weak convergence. Convergence of $R_1$ and $R_2$ in the points $(x, y) \in (0, 1]^2$ for which $H_{\theta_0}(\{y/(x + y)\}) = 0$ follows from Lemma 7.1; see, for instance, [3], Theorem 5.2(iii) (note that by the moment constraints (2.4), $H_\theta/2$ is a probability measure). Since $H_{\theta_0}$ can have at most countably many atoms, $R_1$ and $R_2$ converge in all $(x, y) \in (0, 1]^2$, except for at most countably many rays through the origin. $\qquad\square$

**Proof of Corollary 4.3.** By the continuous mapping theorem, it suffices to show that

$$(\Sigma(\hat{\theta}_n))^{-1/2} D_\varphi(\hat{\theta}_n) \sqrt{k}(\hat{\theta}_n - \theta_0) \xrightarrow{d} N(0, I_p)$$

with $I_p$ being the $p \times p$ identity matrix. By condition (C1) of Theorem 4.2, the map $\theta \mapsto D_\varphi(\theta)$ is continuous at $\theta_0$ so that by the continuous mapping theorem, $D_\varphi(\hat{\theta}_n) \xrightarrow{\mathbb{P}} D_\varphi(\theta_0)$ as $n \to \infty$. Slutsky's lemma and (4.1) yield

$$D_\varphi(\hat{\theta}_n) \sqrt{k}(\hat{\theta}_n - \theta_0) \xrightarrow{d} D_\varphi(\theta_0) D_\varphi(\theta_0)^{-1} \tilde{B} = \tilde{B}$$

as $n \to \infty$. By Lemma 7.2 and the assumption that the map $\theta \mapsto H_\theta$ is weakly continuous, $\Sigma(\hat{\theta}_n)^{-1/2} \xrightarrow{\mathbb{P}} \Sigma(\theta_0)^{-1/2}$. Applying Slutsky's lemma once more concludes the proof. $\qquad\square$

**Proof of Theorem 4.4.** We will show that for the Skorokhod construction introduced before the proof of Theorem 4.2,

$$\left| \iint_{[0,1]^2} (k(\hat{l}_n(x, y) - l(x, y; \hat{\theta}_n))^2 - (B(x, y) - D_{l(x,y;\theta)}(\theta_{\mathcal{H}_0}) D_\varphi(\theta_{\mathcal{H}_0})^{-1} \tilde{B})^2) \, \mathrm{d}x \, \mathrm{d}y \right| \xrightarrow{\mathbb{P}} 0$$

as $n \to \infty$. The left-hand side of the previous expression is less than or equal to

$$\sup_{0 \le x, y \le 1} |\sqrt{k}(\hat{l}_n(x, y) - l(x, y; \hat{\theta}_n)) - B(x, y) + D_{l(x,y;\theta)}(\theta_{\mathcal{H}_0}) D_\varphi(\theta_{\mathcal{H}_0})^{-1} \tilde{B}|$$

$$\times \left( \left| \iint_{[0,1]^2} (\sqrt{k}(\hat{l}_n(x, y) - l(x, y; \theta_{\mathcal{H}_0})) + B(x, y)) \, \mathrm{d}x \, \mathrm{d}y \right| \right.$$

$$\left. + \iint_{[0,1]^2} |\sqrt{k}(l(x, y; \theta_{\mathcal{H}_0}) - l(x, y; \hat{\theta}_n)) - D_{l(x,y;\theta)}(\theta_{\mathcal{H}_0}) D_\varphi(\theta_{\mathcal{H}_0})^{-1} \tilde{B}| \, \mathrm{d}x \, \mathrm{d}y \right)$$

$$=: S(I_1 + I_2).$$

From (7.8) with $g \equiv \mathbf{1}, \mathbf{1} \in \mathbb{R}^p$, we have $I_1 \xrightarrow{\mathbb{P}} 0$. We need to prove that $S = O_{\mathbb{P}}(1)$ and $I_2 = o_{\mathbb{P}}(1)$.

**Proof of $S = O_{\mathbb{P}}(1)$.** We have

$$
\begin{aligned}
S \;\leq\; & \sup_{0 \leq x,y \leq 1} |B(x,y)| + \sup_{0 \leq x,y \leq 1} |\sqrt{k}(\hat{l}_n(x,y) - l(x,y;\theta_{\mathcal{H}_0}))| \\
& + \sup_{0 \leq x,y \leq 1} |\sqrt{k}(l(x,y;\theta_{\mathcal{H}_0}) - l(x,y;\hat{\theta}_n)) + D_{l(x,y;\theta)}(\theta_{\mathcal{H}_0}) D_\varphi(\theta_{\mathcal{H}_0})^{-1}\tilde{B}| \\
=: & \sup_{0 \leq x,y \leq 1} |B(x,y)| + S_1 + S_2.
\end{aligned}
$$

From the definition of process $B$, it follows that $|B(x,y)|$ is almost surely bounded. Furthermore, we have

$$
\begin{aligned}
S_1 = & \sup_{0 \leq x,y \leq 1} |\sqrt{k}(\hat{R}_n^1(x,y) - R(x,y;\theta_{\mathcal{H}_0}))| + o(1) \\
\leq & \sup_{0 \leq x,y \leq 1} |\sqrt{k}(\hat{R}_n^1(x,y) - R_n(S_{1n}(x), S_{2n}(y)))| \\
& + \sup_{0 \leq x,y \leq 1} |\sqrt{k}(R_n(S_{1n}(x), S_{2n}(y)) - R(S_{1n}(x), S_{2n}(y);\theta_{\mathcal{H}_0}))| \\
& + \sup_{0 \leq x,y \leq 1} |\sqrt{k}(R(S_{1n}(x), S_{2n}(y);\theta_{\mathcal{H}_0}) - R(x,y;\theta_{\mathcal{H}_0}))| + o(1)
\end{aligned}
$$

almost surely. In the last part of the proof of Theorem 4.2, we have shown that the third term is almost surely bounded and by the proof of [7], Theorem 2.2, we know that the first two terms are bounded in probability. Let $M$ denote a constant (depending on $\theta_{\mathcal{H}_0}$) bounding the gradient of $\theta \to l(x,y;\theta)$ at $\theta_{\mathcal{H}_0}$ in $(x,y) \in [0,1]^2$. Then, by (4.1),

$$
S_2 \leq M \|\sqrt{k}(\hat{\theta}_n - \theta_{\mathcal{H}_0})\| + M \|D_\varphi(\theta_{\mathcal{H}_0})^{-1}\tilde{B}\| = O_{\mathbb{P}}(1).
$$

**Proof of $I_2 = o_{\mathbb{P}}(1)$.** In Theorem 4.2, we have shown that

$$
T_n := \sqrt{k}(\hat{\theta}_n - \theta_{\mathcal{H}_0}) \xrightarrow{\mathbb{P}} -D_\varphi(\theta_{\mathcal{H}_0})^{-1}\tilde{B} =: N.
$$

By Slutsky's lemma, it is also true that $(T_n, N) \xrightarrow{\mathbb{P}} (N, N)$. By the Skorokhod construction, there exists a probability space, call it $\Omega^*$, which contains both $T_n^*$ and $N^*$, where $(T_n^*, N^*) \stackrel{d}{=} (T_n, N)$ and

$$
(T_n^*, N^*) \xrightarrow{a.s.} (N^*, N^*). \tag{7.9}
$$

Set $\hat{\theta}_n^* := T_n^*/\sqrt{k} + \theta_{\mathcal{H}_0} \stackrel{d}{=} T_n/\sqrt{k} + \theta_{\mathcal{H}_0} = \hat{\theta}_n$. Let $\Omega_0^* \subset \Omega^*$ be a set of probability 1 on which $N^*$ is finite and the convergence in (7.9) holds. We will show that on $\Omega_0^*$,

$$
I_2^* := \iint_{[0,1]^2} X_n^*(x,y)\, \mathrm{d}x\, \mathrm{d}y
$$

$$:= \iint_{[0,1]^2} |\sqrt{k}(l(x,y;\hat{\theta}_n^*) - l(x,y;\theta_{\mathcal{H}_0})) - D_{l(x,y;\theta)}(\theta_{\mathcal{H}_0})N^*| \,\mathrm{d}x\,\mathrm{d}y$$

converges to zero. Since $I_2^* \overset{d}{=} I_2$, the above convergence (namely $I_2^* \overset{a.s.}{\to} 0$) will imply that $I_2 \overset{\mathbb{P}}{\to} 0$. To show that $I_2^*$ converges to zero on $\Omega_0^*$, we will once more apply the dominated convergence theorem. Hereafter, we work on $\Omega_0^*$.

(i) *Pointwise convergence of $X_n^*(x,y)$ to zero.* We have that

$$X_n^*(x,y) \leq |\sqrt{k}(l(x,y;\hat{\theta}_n^*) - l(x,y;\theta_{\mathcal{H}_0}) - D_{l(x,y;\theta)}(\theta_{\mathcal{H}_0})(\hat{\theta}_n^* - \theta_{\mathcal{H}_0}))|$$
$$+ |D_{l(x,y;\theta)}(\theta_{\mathcal{H}_0})(T_n^* - N^*)|.$$

Because of (7.9), differentiability of $\theta \mapsto l(x,y;\theta)$ and continuity of matrix multiplication, the right-hand side of the above inequality converges to zero for all $(x,y) \in [0,1]^2$.

(ii) *Domination of $X_n^*(x,y)$.* Let $M$ be as above. Since the sequence $(T_n^*) = (\sqrt{k}(\hat{\theta}_n^* - \theta_{\mathcal{H}_0}))$ is convergent, and hence bounded, we have

$$\sup_{0 \leq x,y \leq 1} X_n^*(x,y) \leq M\|\sqrt{k}(\hat{\theta}_n^* - \theta_{\mathcal{H}_0})\| + M\|N^*\| = O(1).$$

This concludes the proof of domination and hence the proof of $I_2 \overset{\mathbb{P}}{\to} 0$. $\qquad\square$

**Proof of Lemma 6.1.** Without loss of generality, we can assume that $F$ is itself a parallel elliptical distribution, that is, $(X, Y)$ is given as in (6.1) with $\rho = 0$ in (6.2). Under the assumptions of the lemma and by [18], Theorem 2.3, there exists a function $h: [0, \infty)^2 \to \mathbb{R}$ such that as $t \downarrow 0$ and for all $(x, y) \in [0, \infty)^2$,

$$\frac{t^{-1}\mathbb{P}\{1 - F_1(X) \leq tx, 1 - F_2(Y) \leq ty\} - R(x,y;\nu)}{A(F_2^{\leftarrow}(1-t))} \to h(x,y). \tag{7.10}$$

Moreover, the convergence in (7.10) holds uniformly on $\{(x,y) \in [0,\infty)^2 : x^2 + y^2 = 1\}$ and the function $h$ is bounded on that region; see [18] for an explicit expression of the function $h$.

Condition (6.4) obviously implies that $z \mapsto \mathbb{P}(Z > z)$ is regularly varying at infinity with index $-\nu$. Hence, the same is true for the function $1 - F_2$; see [16]. By [2], Proposition 1.5.7 and Theorem 1.5.12, the function $x \mapsto |A(F_2^{\leftarrow}(1-1/x))|$ is regularly varying at infinity with index $\beta/\nu$. Hence, for every $\alpha < -\beta/\nu$, we have $A(F_2^{\leftarrow}(1-1/x)) = o(x^{-\alpha})$ as $x \to \infty$ or $A(F_2^{\leftarrow}(1-t)) = o(t^\alpha)$ as $t \downarrow 0$. As a consequence, for every $\alpha < -\beta/\nu$, we have, as $t \downarrow 0$,

$$t^{-1}\mathbb{P}\{1 - F_1(X) \leq tx, 1 - F_2(Y) \leq ty\} - R(x,y;\nu) = O(t^\alpha),$$

uniformly on $\{(x,y) \in [0,\infty)^2 : x^2 + y^2 = 1\}$. Uniformity on $\{(x,y) \in [0,\infty)^2 : x+y = 1\}$ now follows as in the proof of [7], Theorem 2.2. $\qquad\square$

# Acknowledgements

The research of Andrea Krajina is supported by an Open Competition grant from the Netherlands Organization for Scientific Research (NWO). The research of Johan Segers is supported by the IAP research network grant no. P6/03 of the Belgian government (Belgian Science Policy) and by a CentER Extramural Research Fellowship from Tilburg University. We are grateful to a referee for providing the reference for the proof of Lemma 7.1.